\documentclass{elsart}
\def\D{\Delta}

\def\a{\alpha}
\def\r{\rho}
\def\t{\tau}
\def\th{\theta}
\def\p{\pi}
\def\l{\ell}

\def\e{\varepsilon}
\def\FF{\mathop{\mathcal F}\nolimits}
\def\suml   {\mathop{\sum}   \limits}
\def\1n{1,\ldots,n}
\def\_#1{\mathop{\hspace{-2pt}^{}_{#1}}}
\def\cdc{,\ldots,}
\def\R{{\rm I\! R}}
\def\Lp {\mathop {L^{\scriptscriptstyle +}}}
\def\lp {\ell^{\scriptscriptstyle +}}
\def\intercal{\mathop {\scriptscriptstyle T}\nolimits}

\begin{document}

\begin{frontmatter}
\title{The Forest Metrics for Graph Vertices$\!\!\!$\thanksref{thank}
}
\thanks[thank]{%
This work was supported by the Russian Foundation for Basic
Research.%
}
\vspace{-2em}

\author{Pavel Chebotarev}$\!$\footnote{Corresponding author. E-mail: {\tt pavel4e@gmail.com, pchv@rambler.ru}}
\and \author{Elena Shamis}
\address{Trapeznikov Institute of Control Sciences of the Russian Academy
of Sciences\\ 65 Profsoyuznaya Street, Moscow 117997, Russia}

\date{April 12, 2002}

\begin{abstract}
We propose a new graph metric and study its properties.
In contrast to the standard distance in connected graphs%
~\cite{BucHar}%
,
it takes into account all paths between vertices.
Formally, it is defined as $d(i,j)=q_{ii}+q_{jj}-q_{ij}-q_{ji}$
\cite{weAiT98}, where $q_{ij}$ is the $(i,j)$-entry of
the {\em relative forest accessibility matrix}
$Q(\varepsilon)=(I+\varepsilon L)^{-1}$, $L$ is the Laplacian matrix of
the (weighted) (multi)graph, and $\varepsilon$ is a positive
parameter. By the matrix-forest theorem, the $(i,j)$-entry of the
relative forest accessibility matrix of a graph provides the specific
number of spanning rooted forests such that $i$ and $j$ belong to the
same tree rooted at $i$. Extremely simple formulas express the
modification of the proposed distance under the basic graph
transformations. We give a topological interpretation of $d(i,j)$ in
terms of the probability of unsuccessful linking $i$ and $j$ in a
model of random links. The properties of this metric are
compared with those of some other graph metrics~\cite{Kle97,weAiT98a}.
An application of this metric is related to clustering procedures such
as centered partition~\cite{Lenart98}. In another procedure, the
relative forest accessibility and the corresponding distance serve to
choose the centers of the clusters and to assign a cluster to each
non-central vertex. Some related geometric representations are
discussed in~\cite{SiSoBi97}. The notion of cumulative weight of connections
between two vertices is proposed. The reasoning involves a reciprocity
principle for weighted multigraphs. Connections between the
resistance distance and the forest distances are established.
\end{abstract}
\end{frontmatter}
\vspace{-2em}

\section{Introduction}

Proximity measures for graph vertices and related algebraic
indices have a wide range of applications. These fall into
information transmission, organic chemistry, crystallography,
integrated circuit design, urban planning, transport networks,
social networks, politology, aggregation of preferences,
epidemiology, etc. (some references are given in~\cite{weAiT98a}).
Let us also mention a few more fields of application, namely,
cluster analysis (see, e.g., \cite{Lenart98}), the theory of
parallel computations \cite{SiSoBi97}, and optics~\cite{Lannes}. In
all these areas, there is a need in proximity measures different
from the classical geodesic distance~\cite{BucHar}.

\setcounter{footnote}{1}
Let $G$ be a weighted multigraph with vertex set $V(G)=\{1\cdc
n\}$ and edge set $E(G)$. Multiple edges are allowed, but loops
are not; $\e_{ij}^p$ is the weight of the $p$th edge between $i$
and $j$; the weights of all edges are strictly
positive.\footnote{In this paper, we consider only such
multigraphs; sometimes we will simply call them graphs.}

In \cite{weAiT97} we proposed the {\em relative forest
accessibility\/} measure of vertex proximity. More exactly, this
is a one-parametric family of indices: if a parameter $\a>0$ is
chosen, then the matrix $Q\_{\a}=(q_{ij}^{\a})$ of vertex
comparative proximities is given by
\begin{equation}
Q\_{\a}=(I+\a L)^{-1},
\label{Q(a)}
\end{equation}
where $I$ is the identity matrix,
$L=L(G)=(\l\_{ij})$ is the Laplacian matrix of the weighted multigraph
$G$ (also termed as the Kirchhoff and the admittance matrix):
\begin{eqnarray}
\label{lij}
\l_{ij}&=&-\suml_{p=1}^{a\_{ij}}\e_{ij}^p,\quad j\neq i,\;\:i,j=\1n,
\\
\l_{ii}&=&-\suml_{j\neq i}\l_{ij},\quad i=\1n,
\label{lii}
\end{eqnarray}
and $a\_{ij}$ is the number of edges with terminal vertices $i$ and
$j$. The matrices $Q\_{\a}$ are doubly stochastic, i.e.,
\begin{eqnarray}
q_{ij}^{\a}&\ge& 0,
\label{nonneg}\\
\sum_{k=1}^nq_{ik}^{\a}&=&\sum_{k=1}^nq_{kj}^{\a}=1,\quad i,j=\1n,
\label{sumone}
\end{eqnarray}
and symmetric. The entry $q_{ij}^{\a}$ can be interpreted as the
relative share of the connections between $i$ and $j$ in the
totality of all $i$'s connections with the vertices of $G$. The
parameter $\a$ determines the proportion of taking into account
long and short routes between vertices.

In \cite{weAiT97} we considered a metric defined by means of the
relative forest accessibilities. The distance between vertices $i$
and $j$ of a (weighted) multigraph was defined as
$q_{ii}^1+q_{jj}^1- q_{ij}^1-q_{ji}^1$. In the present paper, we
study two parametric families of forest metrics of a graph. The
elements of these families are proportional, but the properties of
the families, taken as a whole, are different.
\begin{defn}
\ For a given parameter $\a>0$, the value
\begin{equation}
d_{ij}^{\a}=\frac{1}{2}(q_{ii}^{\a}+q_{jj}^{\a}- q_{ij}^{\a}-q_{ji}^{\a}),
            \quad i,j=\1n
\label{metrad}
\end{equation}
will be called the {\em forest distance between $i$ and $j$}; the
value
\begin{equation}
\r_{ij}^{\a}=\a(q_{ii}^{\a}+q_{jj}^{\a}-q_{ij}^{\a}-q_{ji}^{\a}),
             \quad i,j=\1n
\label{metrar}
\end{equation}
will be called the {\em adjusted forest distance between $i$ and $j$}.
\end{defn}

The fact that these functions are metrics follows from Proposition~1
in~\cite{weAiT97}.

The multiplier $1/2$ in (\ref{metrad}) ensures the equality of
scales for the relative forest accessibilities and the forest
distances (see (\ref{DtauDpi}) and the subsequent remark). If the
vertices $i$ and $j$ belong to the same component of $G$, then
Corollary~1 of~\cite{weAiT98a} implies
$\lim_{\a\to\infty}d_{ij}^{\a}=0$. The asymptotic behavior of the
adjusted forest distances is studied below in
Section~\ref{secresi}. In Section~\ref{seckt}, we demonstrate that
the reversed $\r_{ij}^{\a}$ can be used to measure the cumulative
weight of connections between $i$ and~$j$.

One of the aims of this paper is to present and interpret formulas
for the increments of forest distances and relative forest
accessibilities under the basic multigraph transformations. These
expressions are very simple and testify to the intrinsic nature of
the forest metrics. Also, we will give a stochastic interpretation
of the forest metric and establish some connections between the
forest and resistance metrics.

Until the beginning of Section~\ref{secresi}, $\a$ will be a fixed
parameter, so a simplified notation such as $d\_{ij}$, $\r\_{ij}$,
$q\_{ij}$, and $Q$ will be used instead of the complete notation,
$d_{ij}^{\a}$, $\r_{ij}^{\a}$, $q_{ij}^{\a}$, and $Q_{\a}$,
respectively.

Some properties of the forest metric of a graph have been studied
in \cite{weAiT97} and~\cite{Merr98}. It follows from item~5 of
Proposition~7 in \cite{weAiT97} that if the total weight of all
edges that have terminal vertices $k$ and $t$ gets an increment of
$\D\e\_{kt}>0$, then, provided that there are no other changes,
the increments of the adjusted forest distances in $G$ are
\begin{equation}
\D \r\_{ij}
=-\frac{(\r\_{ik}-\r\_{it}+\r\_{jt}-\r\_{jk})^2}
   {4(\r\_{kt}+1/\D\e_{kt})}, \quad i,j=\1n.
\label{delra}
\end{equation}
Observe that for fixed $\r\_{ik}-\r\_{it}+\r\_{jt}-\r\_{jk}$ and
$\r\_{kt}$, the increments $\D \r\_{ij}$ are independent of $\a$.
Besides that, (\ref{delra}) implies that the forest distances
cannot increase as a result of the addition of a new edge.

The double stochastic property of $Q$ yields the inequalities
\begin{eqnarray}
 d_{ij}&\le& 1,
\label{ine1}\\
\r_{ij}&\le&2\a,\quad i,j=\1n,
\label{ine2}
\end{eqnarray}
and, by proposition~5 from \cite{weAiT97}, the equalities in
(\ref{ine1}) and (\ref{ine2}) are attained if and only if $i$ and
$j$ are two isolated vertices.

Corollary~9 of \cite{Merr97} provides a more accurate upper bound
for $d\_{ij}$:
\begin{equation}
d\_{ij}\le (1+\a a(G))^{-1},
\label{alco}
\end{equation}
where $a(G)$ is  Fiedler's algebraic connectivity of the graph
(which is the second minimal eigenvalue of $L$). Then $(1+\a
a(G))^{-1}$ is, obviously, the second maximal eigenvalue of $Q$.
Upper bounds for the algebraic connectivity of the graph can be
obtained by transforming (\ref{alco}) (see (14) in~\cite{Merr98}).

The {\em diameter\/} of a graph is the greatest distance between
its vertices. This way, (\ref{ine1}) and (\ref{alco}) provide
upper bounds for the forest diameter of a graph. Other bounds for
the forest distances will be given below (equations (\ref{ine3})
and~(\ref{ine4})).

\section{How forest distances change when connections between vertices
strengthen}
\label{seckt}

An expression for the increments of relative forest
accessibilities under basic graph transformations was given in
\cite{weAiT97}, item~1 of Proposition~7 (see also related results
in \cite{Merr98,Merr97,So}). Now we formulate this property in a
more general form.

\begin{defn}
\label{d_odno} \ We say that a weighted multigraph $G'$ {\em
differs from $G$ in a $(k,t)$ edge only\/} if for some
$\D\e\_{kt}=\e\ne0,$ $G'$ can be obtained from $G$ by increasing
the weight of some edge $\e_{kt}^p$ by $\e$, or by adding a new
edge between $k$ and $t$ with weight $\e$ $(\e>0)$, or by removing
an edge between $k$ and $t$ with weight $-\e$ $(\e<0)$.
\end{defn}

While saying that $G'$ differs from $G$ in one edge only, by
$\D\e\_{kt}=\e$ we will denote the weight (or the increment of
weight) in the above definition; primed notation will relate $G'$,
nonprimed notation will relate $G$.

\begin{prop}
\label{prop1} \ Suppose that weighted multigraph $G'$ differs from
$G$ in a $(k,t)$ edge only. Then for every $i,j\in V(G),$
\begin{equation}
 \D q\_{ij}
=q'_{ij}-q\_{ij} =\frac{\a(q\_{ik}-q\_{it})(q\_{jt}-q\_{jk})}
                       {\r\_{kt}+\e^{-1}}.
\label{delq}
\end{equation}
\end{prop}

The results of this section rely on the following lemma, which follows
from (\ref{delra}) after substituting $i=k$ and $j=t$.

\begin{lem}
\label{lem1}
\ Let $G'$ differ from $G$ in a $(k,t)$ edge only. Then
\begin{equation}
\frac{1}{\r'_{kt}}-\frac{1}{\r\_{kt}}=\D\e\_{kt}.
\label{DeltaD}
\end{equation}
\end{lem}

{\bf Proof of Proposition \ref{prop1}} is carried out by the same
argument as the proof of item~1 of Proposition~7
in~\cite{weAiT97}. The only necessary addition is that the
denominator on the right-hand side of (\ref{delq}) cannot vanish
even with a negative $\D\e\_{kt}$.  Indeed, by Lemma~\ref{lem1},
$1/\r\_{kt}+\D\e\_{kt}=1/\r'_{kt}$. As was demonstrated in
\cite{weAiT97}, $\r'_{kt}$ is nonzero (and exists!) at $k\ne t$
for every multigraph $G$, consequently,
$1/\r\_{kt}+\D\e\_{kt}\ne0$, therefore,
$\r\_{kt}+1/\D\e\_{kt}\ne0$.
\bigskip

The connection of an inverse form between the distance and the
increment of weight in (\ref{DeltaD}) requires some explanation.
Note that since $\D\e\_{kt}$ can be arbitrarily large, whereas the
distance between $k$ and $t$ should decrease with any increase of
$\D\e\_{kt}$, a linear connection between the distance and the
$\D\e\_{kt}$ cannot be expected. Taking into account that the
equality $\r'_{kt}=\r\_{kt}$ is inevitable when $\D\e\_{kt}=0$ and
that additivity of distance increments and {\em reciprocity
principle\/} (see below) should also be satisfied, it is difficult
to expect for the distance increments a simpler form
than~(\ref{DeltaD}). It should be also remarked that the inverse
connection between the weight of an edge and its length is the
most natural way of extending the geodesic distance to weighted
graphs (cf.\ \cite{KirNeuSha}).

The distance (\ref{metrar}) for each pair of vertices generally
depends on the graph as a whole, but if $\r\_{kt}$ is fixed, then,
by Lemma~\ref{lem1}, $\r'_{kt}$ is determined by $\D\e\_{kt}$
only. In particular, for any unweighted multigraph, the addition
of a new $(k,t)$ edge yields
\begin{equation}
\frac{1}{\r'_{kt}}-\frac{1}{\r\_{kt}}=1,
\label{DeltaD1}
\end{equation}
which implies
\begin{equation}
\D \r\_{kt}=-\r\_{kt}\r'_{kt}.
\label{lapida}
\end{equation}

The form of (\ref{DeltaD}) and (\ref{DeltaD1}) suggests the
consideration of inverted distances.

\begin{defn}
\label{d_sila} \ Let $i$ and $j\ne i$ be distinct vertices of $G$.
The value $\th\_{ij}=(\r\_{ij})^{-1}-(2\a)^{-1}$ will be called
the {\em cumulative weight of connections between $i$ and $j$ in
$G$}.
\end{defn}

By Lemma~\ref{lem1}, $\th'_{kt}=\th\_{kt}+\D\e\_{kt}$, and
$\th'_{kt}=\th\_{kt}+1$ for unweighted graphs. Successively
applying Lemma~\ref{lem1} to all the edges with terminal vertices
$i$ and $j$ and using (\ref{ine2}) along with its equality
condition, we obtain the following statement, which justifies the
term ``cumulative weight of connections.''

\begin{prop}
\label{prop1'5} \ Suppose that $G$ is an arbitrary weighted
multigraph$,$ $i$ and $j$ are distinct vertices of $G$, and
$\e\_{ij}$ is the total weight of all edges with terminal vertices
$i$ and $j$. Then\\
{\rm 1.} $\th\_{ij}=\th^0_{ij}+\e\_{ij},$
where $\th^0_{ij}$ is the cumulative weight of connections between
$i$ and $j$ in the multigraph resulting from $G$ by the removal of
all edges with terminal vertices $i$ and $j.$\\
{\rm 2.}
$\th\_{ij}\ge\e\_{ij}.$\\
{\rm 3.} $\th\_{ij}=\e\_{ij}$ if and only if $i$ and $j$ are not
connected with the other vertices $($but can be connected to each
other$).$
\end{prop}

Item~2 of Proposition~\ref{prop1'5} and Definition~\ref{d_sila}
provide the following upper bounds for the forest distances:
\begin{eqnarray}
 d_{ij}&\le& (1+2\a\e\_{ij})^{-1},
\label{ine3}\\
\r_{ij}&\le& (\e\_{ij}+(2\a)^{-1})^{-1},\quad i,j=\1n,
\label{ine4}
\end{eqnarray}
where $\e\_{ij}$ is the total weight of the edges with terminal
vertices $i$ and $j$ in~$G$. The equality condition is the same as
in item~3 of Proposition~\ref{prop1'5}.

\section{Comparing the forest distances in two graphs that differ in
one edge}

In this section, we study how a modification on a $(k,t)$ edge alters
the whole profiles of forest distances and relative forest
accessibilities. The alteration of the values
\begin{eqnarray}
\t\_{i(kt)}&=&d\_{ik}-d\_{it},
\label{tau}\\
\p\_{i(kt)}&=&q\_{ik}-q\_{it}, \quad i=\1n
\label{pi}
\end{eqnarray}
will be of major interest here.

In particular, we have
\[
\t\_{t(kt)}=d\_{kt}=-\t\_{k(kt)}.
\]

The simultaneous consideration of forest distances and relative forest
accessibilities rises
a
question: what is the essential difference between them, except
for the fact that a shorter distance usually corresponds to a
greater accessibility? The main difference is in the very {\em
relative\/} nature of the forest accessibility measure: as follows
from Proposition~7 in \cite{weAiT97}, no distance can increase
after the addition of a new edge, whereas the increments of a
relative forest accessibility can have either sign. More
specifically, if $k$ is ``more accessible'' from $i$ than $t$
($q\_{ik}>q\_{it}$), whereas $t$ is ``more accessible'' from $j$
than $k$ ($q\_{jt}>q\_{jk}$), then the addition of a new edge
between $k$ and $t$ increases the accessibility of $j$ from $i$
($q'_{ij}>q\_{ij}$), because the new edge ``widens a road''
between $i$ and~$j$. Inversely, if from both $i$ and $j$, $k$ is
``more accessible'' than $t$ ($q\_{ik}>q\_{it},\,
q\_{jk}>q\_{jt}$), then the addition of a $(k,t)$ edge decreases
the {\em relative\/} forest accessibility of $j$ from $i$
($q'_{ij}<q\_{ij}$). Here, the new edge connects $i$ and $j$ with
$t$ and with the vertices situated ``beyond'' $t$ {\em more
intensively\/} than this edge connects $i$ with $j$ (recall that
the relative forest accessibilities of all vertices from a given
vertex sums to~1). The diagonal entries of $Q$ measure the
``solitariness" of the vertices (which is substantiated by
Theorem~3 in \cite{Merr98}); no one of them can increase when a
new edge is added.

There exists a one-to-one correspondence and a specific duality
\cite{weAiT98} between metrics and $\Sigma$-proximities (one of
which is relative forest accessibility).

Definitions (\ref{metrad}), (\ref{tau}), and (\ref{pi}) imply
simple connections between $d\_{kt}$, $\p\_{i(kt)}$, and
$\t\_{i(kt)}$:
\begin{eqnarray}
2d\_{kt}
&=&\p\_{k(kt)}-\p\_{t(kt)},
\label{distpi}\\
2\t\_{i(kt)}
&=&(\p\_{k(kt)}-\p\_{i(kt)})+(\p\_{t(kt)}-\p\_{i(kt)}),
\label{taupi}\\
\t\_{i(kt)}-\t\_{j(kt)}
&=&\p\_{j(kt)}-\p\_{i(kt)}, \quad i,j,k,t=\1n.
\label{DtauDpi}
\end{eqnarray}

\begin{rem}
Equation (\ref{DtauDpi}) can be interpreted as the equality of
scales of relative forest accessibilities $q\_{ij}$ and forest
distances $d\_{ij}$ (see definitions (\ref{tau}) and (\ref{pi})).
This equality of scales was attained by putting the coefficient
1/2 in~(\ref{metrad}).
\end{rem}

The following statement provides simple expressions for the
increments of all forest distances and accessibilities when one
edge is altered.

\begin{prop}
\ Suppose that $G'$ differs from $G$ in a $(k,t)$ edge only, $\D
d\_{ij}=d'_{ij}-d\_{ij},$ and $\D q\_{ij}=q'_{ij}-q\_{ij}$. Then
the following equalities are true for all $i,j=\1n:$
\begin{eqnarray}
\D q\_{ij}
&=&-\a\e \p\_{i(kt)}\p\_{j(kt)}    \frac{d'_{kt}}{d\_{kt}}
   =   -\a\e \p'_{i(kt)}\p'_{j(kt)}    \frac{d\_{kt}}{d'_{kt}},
\label{Dprox}\\
2\D d\_{ij}
&=&-\a\e(\t\_{i(kt)}-\t\_{j(kt)})^2\frac{d'_{kt}}{d\_{kt}}
   =   -\a\e(\t'_{i(kt)}-\t'_{j(kt)})^2\frac{d\_{kt}}{d'_{kt}}
\nonumber\\
&=&-\a\e(\p\_{i(kt)}-\p\_{j(kt)})^2\frac{d'_{kt}}{d\_{kt}}
   =   -\a\e(\p'_{i(kt)}-\p'_{j(kt)})^2\frac{d\_{kt}}{d'_{kt}},
\label{Ddist}
\end{eqnarray}
where $\e=\D\e\_{kt}$ is the weight of a new edge or the increment
of the weight of an existing edge between $k$ and $t$.
\label{prop2}
\end{prop}

The first equalities in (\ref{Dprox}) and (\ref{Ddist}) are proved
by the application of Lemma~\ref{lem1} to (\ref{delq}) and
(\ref{delra}), respectively. The third equality in (\ref{Ddist})
follows from~(\ref{DtauDpi}). The other connections can be derived
by applying the following {\em reciprocity principle\/}:
\bigskip

{\bf Reciprocity principle.} Suppose that a statement $A$ is true
for every pair of multigraphs $(G,G')$ such that $G'$ differs from
$G$ in a $(k,t)$ edge only. Suppose that statement $A'$ results
from $A$ by replacing all values relating to $G$ with the
corresponding values relating to $G'$ and vice versa and
interchanging $\D\e\_{kt}$ with $\D\e\_{tk}$. Then $A'$ is also
true for every pair of multigraphs $(G,G')$ such that $G'$ differs
from $G$ in a $(k,t)$ edge only.
\bigskip

This principle resembles the oriented duality principle, but its
nature is different. The reciprocity principle follows immediately
from the following obvious fact: if $G'$ differs from $G$ in a
$(k,t)$ edge only, then $G$ also differs from $G'$ in a $(k,t)$
edge only, but with the reversed increment of weight. The
reciprocity principle, in spite of its obviousness, can be used to
obtain some nonobvious corollaries (for example, the equalities in
Proposition~\ref{prop2}).

Let us specify now how the values $\t\_{i(kt)},$ $\p\_{i(kt)}$ change
under the basic transformations of a multigraph.

\begin{prop}
\ Let weighted multigraph $G'$ differ from $G$ in a $(k,t)$ edge
only. Then for all $i=\1n,$ the following equations hold$:$
\begin{eqnarray}
\frac{\t'_{i(kt)}}{d'_{kt}}&=&\frac{\t\_{i(kt)}}{d\_{kt}}
\quad\mbox{or$,$ equivalently$,$}\quad
\frac{d'_{ik}-d'_{it}}{d'_{kt}}=\frac{d\_{ik}-d\_{it}}{d\_{kt}};
\label{proptau}\\
\frac{\p'_{i(kt)}}{d'_{kt}}&=&\frac{\p\_{i(kt)}}{d\_{kt}}
\quad\mbox{or$,$ equivalently$,$}\quad
\frac{q'_{ik}-q'_{it}}{d'_{kt}}=\frac{q\_{ik}-q\_{it}}{d\_{kt}}.
\label{proppi}
\end{eqnarray}
\label{prop3}
\end{prop}

Equation (\ref{proppi}) is a corollary from Proposition~\ref{prop1}
and Lemma~\ref{lem1}, but it can be alternatively obtained from the
equations in (\ref{Dprox}) and (\ref{Ddist}) that are proved by means
of reciprocity principle; (\ref{proptau}) follows, say, from
(\ref{proppi}) and (\ref{taupi}).

According to (\ref{proptau}), the difference of distances
$d\_{ik}-d\_{it}$ reduces as $k$ and $t$ come closer,\footnote{%
This ``coming closer'' is caused by strengthening connections
between $k$ and $t$, whereas the other connections remain the
same.}
so as the proportion in (\ref{proptau}) is preserved, which
makes sense; according to (\ref{proppi}), the same property holds
true for the difference of accessibilities.

Observe also that, by (\ref{proptau})--(\ref{proppi}) and
Lemma~\ref{lem1}, the ratio
\begin{equation}
 \frac{\t'_{i(kt)}}{\t\_{i(kt)}}
=\frac{\p'_{i(kt)}}{\p\_{i(kt)}}
=\frac{d'_{kt}}{d\_{kt}}
=\frac{\r'_{kt}}{\r\_{kt}}
=\frac{1}{ 1+\e \r\_{kt}}
=1-\e \r'_{kt}
\label{proppitau}
\end{equation}
is the same for all $i=\1n$.

\section{An interpretation of the forest metric}

In this section, we give an interpretation of the forest distance
between $i$ and $j$ in terms of spanning rooted forests that connect
$i$ and $j$.

Recall that a rooted forest is an acyclic graph with one vertex marked
as a root in each its component. The components of a rooted forest are
rooted trees.

Every spanning rooted forest of $G$ where $i$ is a root of a tree and
$j$ does not belong to this tree will be called an
{\em unsuccessful connection from $i$ to $j$}. A
{\em   successful connection from $i$ to $j$\/} is a spanning rooted
forest of $G$ where $j$ belongs to a tree rooted at $i$.

By $G\_{\a}$ we denote the weighted multigraph obtained from $G$ by
multiplying the weights of its edges by $\a$.

Consider the following model of choosing a random connection.
\bigskip

{\bf A stochastic model of connecting vertices $i$ and $j$}.
\newline
$1.$ Choose one of the two vertices, $i$ or $j$, with probability $1/2$
each.
\newline
$2.$ Choose a spanning rooted forest of $G\_{\a}$ at random: the
probability of choosing forest $F$ is
\begin{equation}
p(F)=\frac{\e(F)}
          {\suml_{F'\in\FF(G\_{\a})}\e(F')},
\label{proba}
\end{equation}
where $\FF(G\_{\a})$ is the set of all spanning rooted forests of
$G\_{\a}$ and $\e(F')$ is the weight of the forest $F',$ which is
defined as the product of the weights of its edges.
\newline
$3.$ If the forest chosen at stage~2 is an unsuccessful connection from
the vertex chosen at stage~1 to the remaining vertex of the
pair $(i,j)$, then we say that an unsuccessful connection between $i$
and $j$ is chosen.
\bigskip

The following statement is an interpretation of the forest distance
between $i$ and $j$.

\begin{prop}
\ The forest distance $d\_{ij}$ is equal to the probability of choosing
an unsuccessful connection between $i$ and $j$ in the above model of
connecting $i$ and~$j$.
\label{prop4}
\end{prop}

Proposition~\ref{prop4} clarifies the concept of forest metric:
two vertices are close to each other if the probability of
choosing an unsuccessful connection between them is small; the
distance is precisely equal to this probability.

\begin{pf}
\ Observe that, by the matrix-forest theorem~\cite{weAiT97},
$q\_{ii}-q\_{ij}$ is the weighted fraction of unsuccessful
connections from $i$ to $j$. Therefore, by (\ref{proba}),
$q\_{ii}-q\_{ij}$ is the probability of choosing an unsuccessful
connection from $i$ to $j$ on stage~2 of the model. Consequently,
by the total probability formula,
$\frac{1}{2}(q\_{ii}-q\_{ij})+\frac{1}{2}(q\_{jj}-q\_{ji})$ is the
probability of choosing an unsuccessful connection between $i$ and
$j$, and, by (\ref{metrad}), this value is equal to $d\_{ij}$.
\end{pf}

Making use of Theorem~3 in \cite{weAiT98a}, an interpretation of
this kind can be formulated for the resistance distance too.
Connections between the resistance metric and the forest metrics
will be considered in the following section.

\section{Relations between the resistance
metric and the forest metrics}
\label{secresi}

According to item~3 of Proposition~\ref{prop1'5}, the index
$\th\_{ij}$ referred to as a cumulative weight of connections
between $i$ and~$j$ has a rather strange property: if $i$ and $j$
belong to different components of $G$, but at least one of these
vertices is not isolated, then $\th\_{ij}>0$. Which connections
can be associated with such a pair of vertices? The following
analysis of relations between the resistance metric and the forest
metrics clarifies the point.

The resistance metric can be defined as the function that
associates a nonnegative number on the extended real line
$\R\cup\{+\infty\}$ with every pair of vertices $i,j\in V(G)$: if
$i$ and $j$ belong to different components of $G$, then
$\widetilde{\r}\_{ij}=+\infty$; otherwise
\begin{equation}
\widetilde{\r}\_{ij}=\lp_{ii}+\lp_{jj}-\lp_{ij}-\lp_{ji},
\label{rezi}
\end{equation}
where $\lp_{ij}\;,i,j=\1n,$ are the entries of the Moore--Penrose
generalized inverse $\Lp$ of the Laplacian matrix $L$. An
electrical interpretation and some properties of the resistance
metric can be found in \cite{KleRan,Kle97,Bapat99,Patti00}.

Consider the asymptotic behavior of the forest metrics $d^{\a}$ and
$\r^{\a}$. As $a\to 0$, the metric $d^{\a}$ tends to the discrete metric:
\begin{equation}
d^0_{ij}=
\cases{
0, &$j=i,$   \cr
1, &$j\ne i,$\cr
}
\label{disc}
\end{equation}
whereas $\r^{\a}$ tends to the zero function.

Let $V\_i$ be the vertex set of the component of $G$ that contains
vertex~$i$. As has been mentioned in Section~1, if $j\in V\_i$,
then $d_{ij}^{\a}\to0$ as $\a\to\infty$. The behavior of $d^{\a}$
and $\r^{\a}$ as $\a\to\infty$ is completely specified by the
following statement, which follows from Corollary~1 and
Proposition~9 of \cite{weAiT98a}.

\begin{prop}
\label{prop9}
\ The limits $ d_{ij}^{\infty}=\lim_{\a\to\infty} d_{ij}^{\a}$ and
             $\r_{ij}^{\infty}=\lim_{\a\to\infty}\r_{ij}^{\a}$ always
exist $($in the latter case$,$ on the extended real line$),$ and
\begin{eqnarray}
d^{\infty}_{ij}&=&
\cases{
0,                                              &$j\in V\_i,$    \cr
\frac{1}{2}(\frac{1}{|V\_i|}+\frac{1}{|V\_j|}), &$j\not\in V\_i,$\cr
}
\label{broth}\\
\r^{\infty}_{ij}&=&\widetilde{\r}\_{ij},\quad i,j=\1n.
\label{foresi}
\end{eqnarray}
\end{prop}

Note that if $j\not\in V\_i$, then the limit of the cumulative weight
of connections $\th\_{ij}$ as $\a\to\infty$ is~0.

Another relationship between the forest metrics and the resistance
metric is also interesting.

\begin{defn}
\label{alpext}
\ A weighted multigraph $G'$ will be called the {\em $\a$-extension\/}
of the weighted multigraph $G$ if\newline
(i) $V(G')=V(G)\cup\{0\}$, and \newline
(ii) the restriction of $E(G')$ to $V(G)$ coincides with $E(G)$, and
the ratio of the corresponding edge weights in $G'$ and in $G$ is $\a$,
and
\newline
(iii) $E(G')$ contains one edge $(0,i)$ of weight~1 for each
vertex $i\in V(G)$.
\end{defn}

\begin{prop}
\label{prop10} \ For every weighted multigraph $G$ and every
$i,j\in V(G),$ $\widetilde{\r}\_{ij}(G')=2d^{\a}_{ij}(G)$ is
true$,$ where $G'$ is the $\a$-extension of $G,$
$\widetilde{\r}\_{ij}(G')$ is the resistance distance on~$G',$ and
$d^{\a}_{ij}(G)$ is the forest distance on $G$ with
parameter~$\a$.
\end{prop}

\begin{pf}
\ 1. The matrix $Q_{\a}^0$ obtained from $Q\_{\a}=(I+\a L)^{-1}$
by the addition of the zero row and zero column corresponding to
the vertex~0 is a generalized inverse of $L'$, the Laplacian
matrix of~$G'$. Indeed, the required equality
\begin{equation}
\label{obob}
L'Q_{\a}^0L'=L',
\end{equation}
is verified straightforwardly by the multiplication of block
matrices.

2. As has been demonstrated in~\cite{Sharpe65}
(see also~\cite{Styan97}), every generalized inverse $H$ of a doubly
centered matrix $Y$ with singularity~1 can be represented as
\begin{equation}
\label{stya}
H=Y^{\scriptscriptstyle +}+ae^{\intercal}+ea^{\intercal},
\end{equation}
where $Y^{\scriptscriptstyle +}$ is the Moore--Penrose generalized
inverse of $Y$, $a$ and $b$ are some vectors, and
$e=(1,1\cdc 1)^{\intercal}$.

3. Observe that the transformation
\begin{equation}
\label{trans}
h'_{ij}=h\_{ii}+h\_{jj}-h\_{ij}-h\_{ji},
\end{equation}
performed with any matrix $H=(h\_{ij})$ that satisfies
(\ref{stya}) and with $Y^{\scriptscriptstyle +}$ gives the same
result. Applying this observation to the matrix $L'$ which has
singularity~1 and making use of~(\ref{obob}), we obtain the
required statement. Mention that the identity of the results of
transformation (\ref{trans}) for all generalized inverses of the
Laplacian matrix of a connected graph has been noticed
in~\cite{Bapat99}.
\end{pf}

We are now in position to discuss the above-mentioned property of
the cumulative weight index $\th\_{ij}$: this weight can be
positive when $i$ and $j$ belong to different components of~$G$;
on the other hand, its limit value corresponding to the resistance
distance is zero.

The objects modeled by the vertices of the same graph usually have
a specific similarity or a common derivation, which causes their
a~priori relationship. The forest metrics enable one to take such
a relationship into account in the model and to control the
strength of the a~priory connections by the parameter~$\a$: an
increase of $\a$ corresponds to a weakening of a~priory
connections. From this point of view, considering the value
$(\r\_{ij})^{-1}$, instead of
$\th\_{ij}=(\r\_{ij})^{-1}-(2\a)^{-1}$, as the cumulative weight
of connections between $i$ and $j$ makes sense too. In accordance
with Proposition~\ref{prop10}, it is not additional connections
between the objects that represent their a~priori relations, but
connections with the introduced ``hidden source'' vertex. This
way, the input structure of connections between the objects
remains intact.

The introduction of a ``hidden source'' frequently enables one to
obtain simpler proofs for graph-theoretic statements. For
instance, the set of spanning rooted forests of any multigraph $G$
can be put into a one-to-one correspondence with the set of
spanning trees of the 1-extension of~$G$; the corresponding
subgraphs have the same weight. This method makes it possible to
reduce some theorems that involve rooted forests to theorems about
trees.

In some models, varying the parameter $\a$ of $\a$-extension
enables one to specify the relative importance of short and long
connections in a graph. For instance, this is the case for the
linear regression analysis of paired comparisons, where the
addition of a ``hidden source'' corresponds to the turn from the
least squares estimates to the ridge estimates of the object
effects. The matrix of the system of linear equations
corresponding to the ridge estimates is nonsingular as distinct
from the matrix of normal equations. In other words, a
regularization procedure applied to the ill-posed least squares
problem leads to the ridge estimates for the object effects. It
turns out that ridge estimates are more adequate in some cases
than the least squares estimates of the object effects in paired
comparisons~\cite{Cheb94}. These ridge estimates can be
represented via relative forest accessibilities~\cite{Sham94}.
\bigskip

The simple facts and relationships presented in this paper
contribute to the analysis of the resistance and forest metrics
and provide some framework for the comparison of various graph
metrics.

\end{document}